\documentclass[12pt]{article}
\usepackage{amsmath,amssymb}
\usepackage{amsthm}

\def\R{\mathbf{R}}

\def\1{\mathbf{1}}

\def\pa{\partial}
\def\ep{\epsilon}
\def\de{\delta}

\newtheorem{theorem}{Theorem}[section]

\newtheorem{remark}{Remark}

\newcommand{\la}{\lambda}

\begin{document}
\title{Mean field game model of corruption \thanks{Supported by RFFI grant  No 14-06-00326}}
\author{
V. N. Kolokoltsov\thanks{Department of Statistics, University of Warwick,
 Coventry CV4 7AL UK,  Email: v.kolokoltsov@warwick.ac.uk and associate member  of Institute of Informatics Problems, FRC CSC RAS}
 and O. A. Malafeyev\thanks{Fac. of Appl. Math. and Control Processes, St.-Petersburg State Univ., Russia}}
\maketitle

\begin{abstract}
A simple model of corruption that takes into account the effect of the interaction
of a large number of agents by both rational decision making and myopic behavior is developed.
Its stationary version turns out to be a rare example of an exactly solvable model
of mean-field-game type. The results show clearly
how the presence of interaction (including social norms) influences the spread of corruption.
\end{abstract}


\section{Introduction}

Analysis of the spread of corruption in bureaucracy is a well recognized area of the application of game theory,
which attracted attention of many researchers. General surveys can be found in \cite{Aidt}, \cite{Jain}, \cite{LeTsi98}.
In his Prize lecture \cite{Hurw}, L. Hurwicz gives a nice introduction in laymen terms of various problems
arising in an attempt to find out 'who will guard the guardians?' and which mechanisms can be exploited to enforce the legal behavior?
In a series of papers \cite{LaMoMaRa09},\cite{LaMoMaRa08} the authors analyze the dynamic game, where entrepreneurs
have to apply to a set of bureaucrats (in a prescribed order) in order to obtain permission for their
business projects; for an approval the bureaucrats ask for bribes with their amounts being considered as strategies
of the bureaucrats. The existence of an intermediary undertaking the contacts with bureaucrats for a fee may moderate
the outcomes of this game referred to as petty corruption, as each bureaucrat is assumed to ask for a small bribe,
so that the large bureaucratic losses of entrepreneurs occur from a large number of bureaucrats. This is a kind
of extension of the classical ultimatum game, as if an entrepreneur declines to pay the required graft, the game stops.
In the series of works \cite{Vasinbook}, \cite{VasinKarUr10}, \cite{Nik14} the authors develop an hierarchical
model of corruption, where the inspectors of each level audit the inspectors of the previous level and report their finding
to the inspector of the next upper level. For a graft they may choose to make a falsified report. The inspector of the highest level
is assumed to be honest but very costly for the government. The strategy of the government is in the optimal determination
of the audits on various levels with the objective to achieve the minimal level of corruption
with the minimal costs. Paper \cite{Stark14} develops a simple model to get an insight into the problem of when
unifying efforts lead to strength or corruption.
In paper \cite{Mal14} the model of a network corruption game is introduced and
analyzed with the dynamics of corrupted services between the entrepreneurs and corrupted bureaucrats propagating via the chain of  intermediary.
In \cite{NgenZac} the dichotomy between public monitoring and governmental corruptive pressure on the growth of economy
was modeled.  In \cite{LeeSig} an evolutionary model of corruption is developed for
ecosystem management and bio-diversity conservation.
Let us mention also the research on the political aspects of corruption developing around the Acton's dictum that 'power corrupts',
where the elections serve usually as a major tool of public control, see \cite{GiovSeid14} and references therein.
Closely related are the so-called inspection games, see surveys e. g. in
\cite{AvSZ2002}, \cite{ACKvSZ1996}, \cite{KoMabook}, \cite{KoPasWei}.

On the other hand, one of the central trend in the modern theory of games and optimal control
 is related to the analysis of systems with a large number of agents providing a strong link
 with the study of interacting particles in statistical mechanics. Therefore it is natural
 to start applying these methods to the games of corruption, which until recently were mostly studied
 by the classical game-theoretic models with two or three players. The model of corruption with a large number
 of agents interacting in a myopic way (agents try to copy a more successful behavior of their peers)
 was developed in \cite{Ko14} as an example of a general model of pressure and resistance that extends
 the approach of evolutionary games to players interacting in response to a pressure executed
  by a distinguished big player (a principal). In the present paper we consider each player
  of a large group to be a rational optimizer thus bringing the model to the realm of mean-field games.

Mean-field games present a quickly developing area of the game theory.
It was initiated by Lasry-Lions \cite{LL2006} and Huang-Malhame-Caines
\cite{HCM3}, \cite{HCM07}, \cite{Hu10},
see \cite{Ba13}, \cite{Ben13}, \cite{GLL2010}, \cite{Gomsurv}, \cite{Cain14} for recent surveys,
as well as \cite{Car15}, \cite{Car13}, \cite{CarD13}, \cite{Gom14}, \cite{KTY14} and references therein.

New trends concern the theory of mean-field games with a major player \cite{Nour13},
the numeric analysis \cite{Achd13}, the risk-sensitive games \cite{TemBas14},
the games with a discrete state space, see \cite{Gom14a}, \cite{BasRa14} and references therein,
as well as the games and control with a centralized controller of a large pool of agents,
see \cite{GaGaLe} and \cite{Ko12}.

Here we develop a concrete stationary mean-field game model with a finite state space of individual players
describing the distribution of corrupted and honest agents under the pressure of both an
incorruptible governmental representative (often referred to, in literature, as 'benevolent principal',
see e. g. \cite{Aidt})
and the 'social norms' of the society.
This model represents a rare example of being exactly solvable, showing in particular
the non-uniqueness of solutions which is widely discussed in the general mean-field game theory. On the other hand,
this example can be used as a natural toy-model to analyze the link between stationary and dynamic models,
again a nontrivial problem widely discussed in the general theory (though we shall not go in this direction here).
From the point of view of the application to corruption, our contribution is in a systematic study of the interaction
of a large number of (potentially corrupted) agents, each one of them being considered as a rational
optimizer. This mean-field-interaction component of our model can be used to enrich the settings
of the majority of papers cited above.

The paper is organized as follows. In the next two sections we present our model and formulate the main results.
Then we discuss its shortcomings and perspectives. The two final sections contain the proofs.

\section{The model and the objectives of analysis}

An agent is supposed to be in one of the three states: honest $H$, corrupted $C$, reserved $R$, where $R$ is the reserved job
of low salary that an agent receives as a punishment if her corrupted behavior is discovered.

The change between $H$ and $C$ is subject to
the decisions of the agents (though the precise time of the execution of their intent is noisy)
the change from $C$ to $R$ are random with distributions depending on the level of the efforts (say, a budget used)
$b$ of the principal (a government representative) invested in chasing a corrupted behavior,
the change $R$ to $H$ (so-to-say, a new recruitment)
may be possible and is included as a random event with a certain rate.

Let $n_H, n_C, n_R$ denote the numbers of agents in the corresponding states
with $N=n_H+n_C+n_R$ the total number of agents.
By a state of the system we shall mean either the $3$-vector $n=(n_H, n_C, n_R)$
or its normalized version $x=(x_H, x_C, x_R)=n/N$.

 The control parameter $u$ of each player in states $H$ or $C$ may have two values,
 $0$ and $1$, meaning that the player is happy
 with her state ($H$ or $C$) or she prefers to switch one to another; there is no control in the state $R$.
 When the updating decision $1$
 is made, the updating effectively occurs with some rates $\la$. 
The recovery rate, that is the rate of change from $R$ to $H$ (we assume that once recruited the agents start
 by being honest) is a given constant $r$.

Apart from taking a rational decision to swap $H$ and $C$, an honest agent can be pushed to become corruptive by her corruptive
peers, the effect being proportional to the fraction of corrupted agents with certain coefficient $q_{inf}$,
which is analogous to the infection rate in epidemiologic models. On the other hand, the honest agents can contribute
to chasing and punishing corrupted behavior, this effect of a desirable social norm being proportional to the fraction of honest agents
with certain coefficient $q_{soc}$. The presence of the coefficients $q_{inf}$, $q_{soc}$
reflecting the social interaction, makes the dynamics of individual agents
dependent on the distribution of other agents, thus bringing the model to the setting of mean-field games.
It is of our major concern to find out how the presence of interaction influences the spread of corruption.

Thus if all agents use the strategy $u_H, u_C \in \{0,1\}$ and the efforts of the principle is $b$,
the evolution of the state $x$  is clearly given by the ODE
\begin{equation}
 \label{eqmainkineticcorr}
\left\{
 \begin{aligned}
& \dot x_R =(b+q_{soc} x_H) x_C -r x_R, \\
& \dot x_H =r x_R -\la (x_H u_H - x_Cu_C)-q_{inf} x_H x_C, \\
& \dot x_C =-(b+q_{soc} x_H) x_C +\la (x_H u_H - x_Cu_C)+q_{inf} x_H x_C.
\end{aligned}
\right.
 \end{equation}

 It is instructive to see how this ODE can be rigorously deduced from the Markov model of interaction.
 Namely, if all agents use the strategy $u_H, u_C \in \{0,1\}$ and the efforts of the principle is $b$,
 the generator of the Markov evolution on the states $n$ is (where the unchanged values in the arguments of $F$ on the r.h.s
 are omitted)
 \[
 L_NF(n_H, n_C, n_R)=n_C (b+q_{soc}\frac{n_H}{N}) F(n_C-1, n_R+1)+n_R r F(n_R-1, n_H+1)
\]
\[
+n_H (\la u_H +q_{inf}\frac{n_C}{N} ) F(n_H-1, n_C+1) + \la n_C u_C F(n_C-1, n_H+1).
\]
For any $N$, this generator describes a Markov chain on
the finite state space $\{n=(n_H,n_C,n_R): n_H+n_C+n_R=N\}$, where any agent, independently of others,
can be recruited with rate $r$ (if in state $R$) or change from $C$ to $H$ or vice versa if desired
(with rate $\la$), and where the change of the state due to binary interactions are taken into account by the terms
containing $q_{soc}$ and $q_{inf}$.

In terms of $x$ the generator $L_NF$ takes the form
\[
 L_NF(n_H, n_C, n_R)=x_C (b+q_{soc} x_H) F(x-e_C/N+e_R/N)+x_R r F(x-e_R/N+e_H/N)
\]
\begin{equation}
 \label{eqmainlimgencorr0}
+x_H (\la u_H+q_{inf}x_C) F(x-e_H/N+e_C/N) + \la x_C u_C F(x-e_C/N+e_H/N),
\end{equation}
where $\{e_j\}$ is the standard basis in $\R^3$. 
If $F$ is a differentiable function, the generator $L_N$ turns to 
\[
 LF(x)=x_C (b+q_{soc}x_H) \left(\frac{\pa F}{\pa x_R}-\frac{\pa F}{\pa x_C}\right)
 +x_R r \left(\frac{\pa F}{\pa x_H}-\frac{\pa F}{\pa x_R}\right)
 \]
\begin{equation}
 \label{eqmainlimgencorr}
+x_H (\la u_H+q_{inf} x_C) \left(\frac{\pa F}{\pa x_C}-\frac{\pa F}{\pa x_H}\right)
 + \la x_C u_C \left(\frac{\pa F}{\pa x_H}-\frac{\pa F}{\pa x_C}\right)
\end{equation}
in the limit $N\to \infty$.
This is a first order partial differential operator and its characteristics are given by the ODE
 \eqref{eqmainkineticcorr}.

 This Markov model is important not only as a tool to derive \eqref{eqmainkineticcorr},
 but it helps to understand the dynamics of individual players
 (in statistical mechanics terms corresponding to the so-called tagged particles),
 which are central for a mean-field game analysis of agents trying to deviate from the behavior of a crowd.
  Namely, if $x(t)$ and $b(t)$ are given, the dynamics of each individual player is the Markov chain on the $3$ states
with the generator
\begin{equation}
 \label{eqmainlimgenindcorr}
\left\{
\begin{aligned}
& L^{ind}g(R)=r (g(H)-g(R)) \\
& L^{ind}g(H)=(\la u^{ind}_H+q_{inf} x_C) (g(C)-g(H)) \\
& L^{ind}g(C)=\la u^{ind}_C(g(H)-g(C))+(b+q_{soc} x_H) (g(R)-g(C))
\end{aligned}
\right.
\end{equation}
depending on the individual control $u^{ind}\in \{0,1\}$, so that $\dot g=L^{ind}g$ is the Kolmogorov backward
equation of this chain.

Assume that an employed agent receives a wage $w_H$ per unit of time and,
if corrupted, an average payoff $w_C$ (that includes $w_H$ plus some additional illegal reward); she has to pay a fine $f$
when her illegal behavior is discovered; the reserved wage for fired agents is $w_R$.
If the distribution of other payers is $x(t)=(x_R,x_H,x_C)(t)$,
the HJB equation describing the optimal payoff $g=g_t$ (starting at time $t$ with time horizon $T$) of an agent is
\begin{equation}
 \label{eqindHJBcorr}
\left\{
\begin{aligned}
& \dot g(R)+ w_R +r (g(H)-g(R))=0 \\
& \dot g(H)+w_H +\max_u (\la u+q_{inf} x_C) (g(C)-g(H)) =0 \\
& \dot g(C)+w_C -(b+q_{soc}x_H)f  +\max_u (\la u (g(H)-g(C))+(b+q_{soc} x_H) (g(R)-g(C))=0.
\end{aligned}
\right.
 \end{equation}

Therefore, starting with some control
\[
u^{com}(t)=(u^{com}_C(t), u^{com}_H(t)),
\]
used by all players, we can find the dynamics $x(t)$ from equation \eqref{eqmainkineticcorr}
(with $u^{com}$ used for $u$). Then each individual should solve
the Markov control problem \eqref{eqindHJBcorr}
thus finding the individually optimal strategy
\[
u^{ind}(t)=(u^{ind}_C(t), u^{ind}_H(t)).
\]

 The {\it basic MFG consistency} equation can now be explicitly written as
 \begin{equation}
 \label{eqMFGconscorr}
 u^{ind}(t)=u^{com}(t).
 \end{equation}

Instead of analyzing this rather complicated dynamic problem, we shall look for a simpler and practically more
relevant problem of consistent stationary strategies.

There are two standard stationary problems arising from
HJB \eqref{eqindHJBcorr}, one being the search for the average payoff
\[
g =\lim_{T\to \infty} \frac{1}{T}\int_0^T g_t \, dt
\]
for long period games, and another the search for discounted optimal payoff.
The first is governed by the solutions of HJB of the form $(T-t)\mu +g$, linear in $t$ (with $\mu$ describing the optimal average payoff),
so that $g$ satisfies the stationary HJB equation:

\begin{equation}
 \label{eqindstHJBcorr}
\left\{\begin{aligned}
& w_R +r (g(H)-g(R))=\mu \\
& w_H +\max_u (\la u+q_{inf} x_C) (g(C)-g(H)) =\mu \\
& w_C -(b+q_{soc} x_H) f +\max_u (\la u (g(H)-g(C))+(b+q_{soc} x_H) (g(R)-g(C))=\mu,
\end{aligned}
\right.
 \end{equation}
 and the discounted optimal payoff (with the discounting coefficient $\de$) satisfies the stationary HJB

 \begin{equation}
 \label{eqindstHJBcorrdis}
\left\{
\begin{aligned}
& w_R +r (g(H)-g(R))=\de g (R) \\
& w_H +\max_u (\la u+q_{inf} x_C) (g(C)-g(H)) =\de g(H) \\
& w_C -(b+q_{soc} x_H) f +\max_u (\la u (g(H)-g(C))+(b+q_{soc} x_H) (g(R)-g(C))=\de g(C).
\end{aligned}
\right.
 \end{equation}

The analysis of these two settings is mostly analogous (as they are in some sense equivalent, see e. g. \cite{Ross}).
 We shall concentrate on the first one.

For a fixed $b$, the {\it stationary MFG consistency} problem is in finding $(x,u_C,u_H)=(x,u_C(x),u_H(x))$, where $x$ is the stationary point of
evolution \eqref{eqmainkineticcorr}, that is
\begin{equation}
 \label{eqindstHJBcorrequi}
 \left\{
\begin{aligned}
& (b+q_{soc} x_H) x_C -r x_R =0 \\
& r x_R -\la (x_H u_H(x) - x_Cu_C(x))-q_{inf} x_H x_C=0 \\
& -(b+q_{soc} x_H) x_C  +\la (x_H u_H(x) - x_Cu_C(x))+q_{inf} x_H x_C=0,
\end{aligned}
\right.
\end{equation}
where $u_C(x), u_H(x)$ give maximum in the solution to \eqref{eqindstHJBcorr}. Thus $x$ is a fixed point of the
limiting dynamics of the distribution of large number of agents such that the corresponding
stationary control is individually optimal subject to this distribution.

\begin{remark}
Notice that our stationary MFG consistency is close to the concept of the Wardrop equilibria,
see e. g. \cite{HauMac85}, but is quite different nevertheless.
\end{remark}

Fixed points can practically model a stationary behavior only if they are stable. Thus we are interested
in stable solutions $(x,u_C,u_H)=(x,u_C(x),u_H(x))$ to the stationary MFG consistency problem, where
a solution is stable if the corresponding stationary distribution $x=(x_R,x_H,x_C)$ is a stable equilibrium to
\eqref{eqmainkineticcorr} (with $u_C,u_H$ fixed by this solution).
As mentioned above, our major concern is to find out how the presence of interaction
(specified by the coefficients $q_{soc}, q_{inf}$) affects the stable equilibria.

\section{Results}

Our first result describe explicitly all solutions to the stationary MFG consistency problem stated above
and the second result deals with the stability of these solutions.

We shall say that in a solution to the stationary MFG consistency problem  the optimal individual behavior is corruption
if $u_C=0, u_H=1$: if you are corrupt stay corrupt, and if you are honest,
start corrupted behavior as soon as possible; the optimal individual behavior is honesty if
$u_C=1, u_H=0$: if you are honest stay honest,
if you are involved in corruption try to clean yourself from corruption as soon as possible.

The basic assumptions on our coefficients are
\begin{equation}
 \label{eqcoefcorrassum}
\la>0, r> 0, b>0, \quad f\ge 0, q_{soc} \ge 0, q_{inf}\ge 0, \quad w_C > w_H > w_R \ge 0.
\end{equation}
The key parameter for our model turns out to be the quantity
\begin{equation}
 \label{eqcoefcorrleyparam}
\bar x =\frac{1}{q_{soc}}\left[ \frac{r(w_C-w_H)}{w_H-w_R+rf}-b\right]
\end{equation}
(which can take values $\pm \infty$ if $q_{soc}=0$).

\begin{theorem}
\label{thsolvestatMFGcorr}
Assume \eqref{eqcoefcorrassum}.

(i) If $\bar x >1$, then there exists a unique solution $x^*=(x^*_R, x^*_C, x^*_H)$ to the stationary MFG problem \eqref{eqindstHJBcorrequi},
\eqref{eqindstHJBcorr}, where
\begin{equation}
 \label{eq1thsolvestatMFGcorr}
x_C^*=\frac{(1-x_H^*)r}{r+b+q_{soc} x_H^*}
\end{equation}
and $x_H^*$ is the unique solution on the interval $(0,1)$ of the quadratic equation $Q(x_H)=0$, where
\begin{equation}
 \label{eq1athsolvestatMFGcorr}
Q(x_H)=[(r+\la) q_{soc} -r q_{inf}] x_H^2 +[r(q_{inf}-q_{soc})+\la r+\la b +r b]x_H -rb.
\end{equation}
Under this solution the optimal individual behavior is corruption: $u_C=0, u_H=1$.

(ii) If  $\bar x <1$, there may be 1,2 or 3 solutions to the stationary MFG problem \eqref{eqindstHJBcorrequi},
\eqref{eqindstHJBcorr}. Namely,
the point $x_H=1, x_C=x_R=0$ is always a solution,
under which the optimal individual behavior is being honest: $u_C=1, u_H=0$.

Moreover, if
\begin{equation}
 \label{eq2thsolvestatMFGcorr}
\max(\bar x,0) \le \frac{b+\la}{q_{inf} -q_{soc}} <1,
\end{equation}
then there is another solution with the optimal individual behavior being honest, that is $u_C=1, u_H=0$:
\begin{equation}
 \label{eq3thsolvestatMFGcorr}
x^{**}_H= \frac{b+\la}{q_{inf} -q_{soc}}, \quad x_C^{**}=\frac{r(q_{inf}-q_{soc}-b-\la)}{(r+b)q_{inf}+(\la -r) q_{soc}}.
\end{equation}

Finally, if
\begin{equation}
 \label{eq4thsolvestatMFGcorr}
\bar x > 0, \quad Q(\bar x)\ge 0,
\end{equation}
there is a solution with the corruptive optimal behavior of the same structure as in (i), that is,
with $x_H^*$ being the unique solution to $Q(x_H)=0$ on $(0,\bar x]$ and $x_C^*$ given by \eqref{eq1thsolvestatMFGcorr}.

\end{theorem}

\begin{remark}
As seen by inspection, $Q[(b+\la)/(q_{inf}-q_{soc})]>0$ (if $q_{inf}-q_{soc}>0$), so that for $\bar x$ slightly
less than $x_H^{**}=(b+\la)/(q_{inf}-q_{soc})$
one has also $Q(\bar x)>0$, in which case one really has three points of equilibria given by $x_H^*, x_H^{**}, x_H=1$
with $0< x^* <\bar x < x^{**} <1$.
\end{remark}

\begin{remark}
In case of the stationary problem arising from the discounting payoff, that is from equation \eqref{eqindstHJBcorrdis},
the role of the classifying parameter $\bar x$ from \eqref{eqcoefcorrleyparam} is played by the quantity

\begin{equation}
 \label{eqcoefcorrleyparamdisc}
\bar x =\frac{1}{q_{soc}}\left[ \frac{(r+\de)(w_C-w_H)}{w_H-w_R+(r+\de)f}-b\right].
\end{equation}
\end{remark}

\begin{theorem}
\label{thsolvestatMFGcorrstable}
Assume \eqref{eqcoefcorrassum}.

(i) The solution $x^*=(x^*_R, x^*_C, x^*_H)$ (given by Theorem \ref{thsolvestatMFGcorr})
with individually optimal behavior being corruption is stable if
\begin{equation}
 \label{eq1thsolvestatMFGcorrstable}
-\frac{\la q_{soc}}{r}\le q_{soc}-q_{inf} \le \frac{rq_{inf} +(r+b)(br +r\la +b\la)}{r^2}.
\end{equation}

(ii) Suppose  $\bar x <1$.
If $q_{inf}-q_{soc}\le 0$ or
\[
 q_{inf}-q_{soc}> 0, \quad \frac{b+\la}{q_{inf} -q_{soc}} >1,
 \]
then $x_H=1$ is the unique stationary MFG solution with individually optimal strategy being honest;
and this solution is stable.
If \eqref{eq2thsolvestatMFGcorr} holds, there are two stationary MFG solution with individually optimal strategy being honest,
one with $x_H=1$ and another with $x_H=x_H^{**}$ given by \eqref{eq3thsolvestatMFGcorr}; the first solution is unstable
and the second is stable.

\end{theorem}

We are not presenting necessary and sufficient condition for the stability of solutions with optimally corrupted behavior.
Condition \eqref{eq1thsolvestatMFGcorrstable} is only sufficient, but it covers a reasonable range of parameters
where the 'epidemic' spread of corruption and social cleaning are effects of comparable order.

As a trivial consequence of our theorems we can conclude that in the absence of interaction,
that is for $q_{inf}=q_{soc}=0$, the corruption is individually optimal if
\begin{equation}
 \label{eq1nointer}
w_C-w_R \ge bf +(w_H-w_R) (1+\frac{b}{r})
\end{equation}
and honesty is individually optimal otherwise (which is of course a reformulation of 
the standard result for a basic model of corruption, see e. g. \cite{Aidt}).
In the first case the unique equilibrium is
\begin{equation}
 \label{eq2nointer}
x_H^*=\frac{rb}{\la r +\la b +rb}, \quad x_C^*=\frac{r(1-x_H^*)}{r+b},
\end{equation}
and in the second case the unique equilibrium is $x_H=1$. Both are stable.

\section{Discussion}

The results above show clearly
how the presence of interaction (including social norms) influences the spread of corruption.
When $q_{inf}=q_{soc}=0$, one has one equilibrium that corresponds to corrupted or honest behavior
depending on a certain relation \eqref{eq1nointer} between the parameters of the game.
If social norms or 'epidemic' myopic behavior are allowed in the model, which is quite natural
for a realistic process, the situation becomes much more complicated.
In particular, in a certain range of parameters, one has two stable equilibria,
one corresponding to an optimally honest and another to an optimally corrupted behavior.
This means in particular that similar strategies of a principal (defined by the choice of parameters $b,f,w_H$)
can lead to quite different outcomes depending on the initial distributions of
honest-corrupted agents or even on the small random fluctuations in the process of evolution.

The coefficients $b$ and $f$ enter exogenously in our system and can be used as tools
for shifting the (precalculated) stable equilibria in the desired direction. These coefficients
are not chosen strategically, which is an appropriate assumption for situations when the principal may have only poor
information about the overall distribution of states of the agents.
It is of course natural to extend the model by treating the principal as a strategic optimizer who
chooses $b$ (or even can choose $f$) in each state to optimize certain payoff. This would place the model
in the group of MFG models with a major player, which is actively studied in the current literature.

Classifying agents as corrupted and honest only is a strong simplification of reality. In the spirit of 
\cite{Nik14} and \cite{Ko14}
it is natural to consider the hierarchy $i=1, \cdots , n$ of the possible positions of agents in a bureaucratic staircase
with both basic wages $w_H^i$ and the illegal payoff $w_C^i$ in the corresponding states $H_i$ and $C_i$
increasing with $i$. Once a corruptive behavior of an agent from
state $C_i$ is detected, she is supposed to be downgraded to the reserved state $R=H_0$, and the upgrading from
$i$ to $i+1$ can be modeled as a random event with a given rate. This {\it multi-layer model of corruption}
could bring insights on the spread of corruption among the representatives of the different levels of power.

Theoretically, the main questions left open by our analysis are the precise link between the stationary and dynamic
 MFG solutions and the precise statement of the law o large numbers. Namely,
 (i) Can we solve the dynamic MFG consistency problem \eqref{eqMFGconscorr} and whether its solutions will approach
 the solutions of the stationary problems described  by our Theorems?
 (ii) Considering a stochastic game of $N$ players in the Markov model where each player evolves according
 to \eqref{eqmainlimgenindcorr} with chosen control $u_C,u_H$ and the distribution $x_t$ reflects
 the aggregated distribution so obtained, do our stationary MFG solutions represent approximate Nash equilibria
 to this game? This latter question is an MFG version of the well known problem of evolutionary game theory
 about the correspondence between the results of taking limits $N\to \infty$ and $t \to \infty$
 in a different order, where rather deep results were obtained, see e. g. \cite{BinSam} and references therein.

\section{Proof of Theorem \ref{thsolvestatMFGcorr}}

Clearly solutions to \eqref{eqindstHJBcorr} are defined up to an additive constant.
Thus we can and will assume that $g(R)=0$. Moreover, we can reduce the analysis to the case
$w_R=0$ by subtracting it from all equations of \eqref{eqindstHJBcorr} and thus shifting by $w_R$ the values
$w_H,w_C, \mu$. Under these simplifications, the first equation to \eqref{eqindstHJBcorr} is $\mu=rg(H)$,
so that \eqref{eqindstHJBcorr} becomes the system
\begin{equation}
 \label{eqindstHJBcorrsim1}
\left\{
\begin{aligned}
& w_H +\la \max(g(C)-g(H),0) +q_{inf} x_C (g(C)-g(H)) =rg(H) \\
& w_C -(b+q_{soc} x_H) f +\la \max(g(H)-g(C),0)-(b+q_{soc} x_H) g(C)=rg(H)
\end{aligned}
\right.
 \end{equation}
 for the pair $(g(H), g(C))$ with $\mu=rg(H)$.

 Assuming $g(C)\ge g(H)$, that is $u_C=0, u_H=1$, so that the corruptive behavior is optimal,
 system \eqref{eqindstHJBcorrsim1} turns to
\begin{equation}
 \label{eqindstHJBcorrsimCop}
\left\{
\begin{aligned}
& w_H +\la (g(C)-g(H)) +q_{inf} x_C (g(C)-g(H)) =rg(H) \\
& w_C -(b+q_{soc} x_H) f -(b+q_{soc} x_H) g(C)=r g(H).
\end{aligned}
\right.
 \end{equation}

Solving this system of two linear equations we get
\[
g(C)=\frac{(r+\la +q_{inf}x_C)[w_C-(b+q_{soc} x_H) f] -rw_H}{r(\la +q_{inf}x_C+b+q_{soc}x_H)+(\la +q_{inf}x_C)(b+q_{soc}x_H)},
\]
\[
g(H)=\frac{(\la +q_{inf}x_C)[w_C-(b+q_{soc} x_H) f]+(b+q_{soc}x_H) w_H}{r(\la +q_{inf}x_C+b+q_{soc}x_H)+(\la +q_{inf}x_C)(b+q_{soc}x_H)},
\]
so that $g(C)\ge g(H)$ is equivalent to
\[
w_C-(b+q_{soc} x_H) f \ge w_H \left(1+\frac{b+q_{soc}x_H}{r}\right),
 \]
 or, in other words,
\begin{equation}
 \label{eqindstHJBcorrsimCopcond1}
x_H \le \frac{1}{q_{soc}}\left[ \frac{r(w_C-w_H)}{w_H+rf}-b\right],
 \end{equation}
 which by restoring $w_R$ (shifting $w_C,w_H$ by $w_R$) gives
\begin{equation}
 \label{eqindstHJBcorrsimCopcond1a}
x_H \le \bar x=\frac{1}{q_{soc}}\left[ \frac{r(w_C-w_H)}{w_H-w_R+rf}-b\right].
 \end{equation}
Since $x_H\in (0,1)$, this is automatically satisfied if $\bar x >1$, that is under the assumption of (i).
On the other hand, it definitely cannot hold if $\bar x <0$.

Assuming $g(C)\le g(H)$, that is $u_C=1, u_H=0$, so that the honest behavior is optimal,
 system \eqref{eqindstHJBcorrsim1} turns to
\begin{equation}
 \label{eqindstHJBcorrsimHop}
\left\{
\begin{aligned}
& w_H +q_{inf} x_C (g(C)-g(H)) =rg(H) \\
& w_C-(b+q_{soc} x_H) f  +\la (g(H)-g(C))-(b+q_{soc} x_H) g(C)=rg(H).
\end{aligned}
\right.
 \end{equation}

Solving this system of two linear equations we get
\[
g(C)=\frac{(r+q_{inf}x_C)[w_C-(b+q_{soc} x_H) f]+(\la -r)w_H}{r(\la +q_{inf}x_C+b+q_{soc}x_H)+q_{inf}x_C(b+q_{soc}x_H)}
\]
\[
g(H)=\frac{q_{inf}x_C [w_C-(b+q_{soc} x_H) f]+(\la + b+q_{soc}x_H) w_H}{r(\la +q_{inf}x_C+b+q_{soc}x_H)+q_{inf}x_C(b+q_{soc}x_H)}
\]
so that $g(C)\le g(H)$ is equivalent to the inverse of condition \eqref{eqindstHJBcorrsimCopcond1}.

If $g(C)\ge g(H)$, that is $u_C=0, u_H=1$, the fixed point equation \eqref{eqindstHJBcorrequi} becomes

\begin{equation}
 \label{eqindstHJBcorrequia}
 \left\{
\begin{aligned}
& (b+q_{soc} x_H) x_C -r x_R =0 \\
& r x_R -\la x_H -q_{inf} x_H x_C=0 \\
& -(b+q_{soc} x_H) x_C +\la x_H +q_{inf} x_H x_C=0.
\end{aligned}
\right.
\end{equation}

Since $x_R=1-x_H-x_C$, the third equation is a consequence of the first two equations, which yields the system
\begin{equation}
 \label{eqindstHJBcorrequi1}
\begin{aligned}
& (b+q_{soc} x_H) x_C -r(1-x_H-x_C)=0 \\
& r (1-x_H-x_C)-\la x_H -q_{inf} x_H x_C=0.
\end{aligned}
\end{equation}
From the first equation we have
\begin{equation}
 \label{eqindstHJBcorrequi1a}
x_C=\frac{(1-x_H)r}{r+b+q_{soc} x_H}.
\end{equation}
From this it is seen that if $x_H\in (0,1)$ (as it should be), then also $x_C\in (0,1)$ and
\[
x_C+x_H=\frac{r+x_H(b+q_{soc} x_H)}{r+b+q_{soc}x_H} \in (0,1).
\]
Plugging $x_C$ in the second equation of \eqref{eqindstHJBcorrequi1} we find for $x_H$
the quadratic equation $Q(x_H)=0$ with $Q$ given by \eqref{eq1athsolvestatMFGcorr}.

Since $Q(0)<0$ and $Q(1)>0$, the equation $Q(x_H)=0$ has exactly one positive root $x_H^*\in (0,1)$.
Hence $x_H^*$ satisfies \eqref{eqindstHJBcorrsimCopcond1} if and only if either $\bar x>1$
(that is we are under the assumption of (i)) or
if \eqref{eq4thsolvestatMFGcorr} holds proving the last statement of (ii).

If $g(C)\le g(H)$, that is $u_C=1, u_H=0$, the fixed point equation \eqref{eqindstHJBcorrequi} becomes

\begin{equation}
 \label{eqindstHJBcorrequiHop}
\left\{
\begin{aligned}
& (b+q_{soc} x_H) x_C -x_R r=0 \\
& x_R r+\la x_C -q_{inf} x_H x_C=0 \\
& -x_C (b+q_{soc} x_H) -\la x_C +q_{inf} x_H x_C=0.
\end{aligned}
\right.
\end{equation}

Again here $x_R=1-x_H-x_C$ and the third equation is a consequence of the first two equations, which yields the system
\begin{equation}
 \label{eqindstHJBcorrequi1Hop}
\left\{
\begin{aligned}
& (b+q_{soc} x_H) x_C -r(1-x_H-x_C)=0 \\
& r (1-x_H-x_C)+\la x_C -q_{inf} x_H x_C=0.
\end{aligned}
\right.
\end{equation}
From the first equation we again get \eqref{eqindstHJBcorrequi1a}.
Plugging this $x_C$ in the second equation of \eqref{eqindstHJBcorrequi1} we find the equation
\[
r(1-x_H)=(r-\la +q_{inf}x_H)\frac{(1-x_H)r}{r+b+q_{inf} x_H},
\]
with two explicit solutions yielding the first and the second statements of (ii).

\section{Proof of Theorem \ref{thsolvestatMFGcorrstable}}

(i) When individually optimal behavior is to be corrupted, that is $u_C=0, u_H=1$,
system \eqref{eqmainkineticcorr} written in terms of $(x_H,x_C)$ becomes

\begin{equation}
 \label{eqmainkineticcorrCop}
 \left\{
 \begin{aligned}
& \dot x_H =(1-x_H-x_C) r-\la x_H-q_{inf} x_H x_C, \\
& \dot x_C =-x_C (b+q_{soc} x_H) +\la x_H +q_{inf} x_H x_C.
\end{aligned}
 \right.
 \end{equation}

 Written in terms of $y=x_H-x_H^*, z=x_C-x_C^*$ it takes the form
\begin{equation}
 \label{eqmainkineticcorrCop1}
 \left\{
 \begin{aligned}
& \dot y =-y( r+\la +q_{inf}x_C^*)-z (r+q_{inf}x_H^*)-q_{inf}yz, \\
& \dot z =y [\la+(q_{inf}-q_{soc}) x_C^*] +z[x_H^*(q_{inf}-q_{soc})-b]z +(q_{inf}-q_{soc})yz.
\end{aligned}
 \right.
 \end{equation}

The condition of stability is the requirement that both eigenvalues of the linear approximation around the fixed point
have real negative parts, or equivalently that the trace of the linear approximation is negative and the determinant is positive:
\begin{equation}
 \label{eqstabilcorrCop1}
 \begin{aligned}
& x_H^*(q_{inf}-q_{soc})-b -r-\la -q_{inf}x_C^* <0, \\
& \la (r+q_{soc} x_H^* +b)-rx_H^*(q_{inf}-q_{soc})+br +x_C^* [r(q_{inf}-q_{soc})+q_{inf}b]>0
\end{aligned}
 \end{equation}
 (note that the quadratic terms in $x_C,x_H$ cancel in the second inequality).
By \eqref{eq1thsolvestatMFGcorr} this rewrites in terms of $x_H^*$ as
 \[
 \begin{aligned}
& [x_H^*(q_{inf}-q_{soc})-b -r-\la] (r+b+q_{soc} x_H^*)-q_{inf}r(1-x_H^*)<0, \\
& [\la (r+q_{soc} x_H^* +b)-rx_H^*(q_{inf}-q_{soc})+br](r+b+q_{soc}x_H^*) +r(1-x_H^*)[r(q_{inf}-q_{soc})+q_{inf}b]>0
\end{aligned}
 \]
 or in a more concise form as
 \begin{equation}
 \label{eqstabilcorrCop2}
 \begin{aligned}
& (x_H^*)^2(q_{inf}-q_{soc})q_{soc} +x_H^*[(q_{inf}-q_{soc})(2r+b)-q_{soc} (b+\la)]-(r+b)(r+b+\la)-rq_{inf}<0, \\
& (x_H^*)^2q_{soc}[(q_{inf}-q_{soc})r-\la q_{soc}]+ 2x_H^*(r+b)[r(q_{inf}-q_{soc})r-\la q_{soc}] \\
& \quad \quad \quad -r^2(q_{inf}-q_{soc})-rbq_{inf}-(r+b)(br+r\la+b\la)<0.
\end{aligned}
 \end{equation}

Let
\[
0 \le q_{soc}-q_{inf} \le \frac{rq_{inf} +(r+b)(br +r\la +b\la)}{r^2}.
\]
Then it is seen directly that both inequalities in \eqref{eqstabilcorrCop2} hold trivially for any positive $x_H$.

Assume now that
\[
0< r(q_{inf}-q_{soc})\le \la q_{soc}.
\]
Then the second condition in \eqref{eqstabilcorrCop2} again holds trivially for any positive $x_H$.
Moreover, it follows from $Q(x_H^*)=0$ that
\[
x_H^* \le \frac{rb}{r (q_{inf}-q_{soc})+\la r +\la b +rb} \le \tilde x=\frac{b}{q_{inf}-q_{soc}}.
\]
Now the left hand side of the first inequality of \eqref{eqstabilcorrCop2} evaluated at $\tilde x$ is negative,
because it equals
\[
-\frac{bq_{soc} \la}{q_{inf}-q_{soc}}-r^2-\la (r+b)-r q_{inf},
\]
and it is also negative when evaluated at $x_H^* \le \tilde x$.

(ii)
When individually optimal behavior is to be honest,
that is $u_C=1, u_H=0$, system \eqref{eqmainkineticcorr} written in terms of $(x_H,x_C)$ becomes

\begin{equation}
 \label{eqmainkineticcorrHop}
 \begin{aligned}
& \dot x_H =(1-x_H-x_C) r+\la x_C-q_{inf} x_H x_C, \\
& \dot x_C =-x_C (b+q_{soc} x_H) -\la x_C +q_{inf} x_H x_C.
\end{aligned}
 \end{equation}
 To analyze the stability of the fixed point $x_H=1, x_C=0$ we write it in terms of $x_C$ and $y=1-x_H$ as
 \[
 \begin{aligned}
& \dot y =-ry +x_C(r-\la +q_{inf}) -q_{inf} y x_C, \\
& \dot x_C =x_C (q_{inf} -q_{soc} -\la -b) -y x_C(q_{inf}-q_{soc}).
\end{aligned}
\]
According to the linear approximation, the fixed point $y=0,x_C=0$ of this system is stable if
$q_{inf} -q_{soc} -\la -b<0$ proving the first statement in (ii).

Assume \eqref{eq2thsolvestatMFGcorr} holds. To analyze the stability of the fixed point $x_H^{**}$ we write
system \eqref{eqmainkineticcorrHop} in terms of the variables
\[
y=x_H-x_H^{**}=x_H-\frac{b+\la}{q_{inf} -q_{soc}},
\quad
z=x_C-x_C^{**}=x_C-\frac{r(q_{inf}-q_{soc}-b-\la)}{(r+b)q_{inf}+(\la -r) q_{soc}},
\]
which is
\[
\begin{aligned}
& \dot y =-y \frac{r[(r+q_{inf})(q_{inf}-q_{soc})+\la q_{soc}]}{(r+b)q_{inf} +(\la -r) q_{soc}}
-z \frac{(r+b)q_{inf} +(\la -r) q_{soc}}{q_{inf}-q_{soc}}-q_{inf} yz, \\
&\dot z =y \frac{r(q_{inf}-q_{soc}-b-\la)(q_{inf}-q_{soc})}{(r+b)q_{inf} +(\la -r) q_{soc}} +yz.
\end{aligned}
\]

The characteristic equation of the matrix of linear approximation is seen to be
\[
\xi^2+ \frac{r[(r+q_{inf})(q_{inf}-q_{soc})+\la q_{soc}]}{(r+b)q_{inf} +(\la -r) q_{soc}}\xi +r(q_{inf}-q_{soc}-b-\la)=0.
\]
Under \eqref{eq2thsolvestatMFGcorr} both the free term and the coefficient at $\xi$ are positive. Hence
 both roots have negative real parts implying stability.

\end{document}